\documentclass[smallextended,numbook,runningheads]{svjour3}     % onecolumn (second format)
\usepackage[svgnames]{xcolor}
\usepackage{amsmath,amssymb,graphicx,hyperref}
\usepackage{booktabs}
\usepackage{mathptmx}      % use Times fonts if available on your TeX system

\def\makeheadbox{\relax}
\RequirePackage{tikz}

\usetikzlibrary{arrows,trees,positioning,calc,matrix}

\usepackage{etrees}
\newcommand{\tr}{\triangleright}
\newcommand{\RR}{\mathbb R}
\newcommand{\F}{{\mathcal F}}

\usetikzlibrary{matrix}
\tikzstyle{commdiag}=[matrix of math nodes, row sep=3em, column sep=2.5em, text height=1.5ex, text depth=0.25ex,ampersand replacement=\&]
\tikzset{>=stealth}
\colorlet{emphcolor}{FireBrick}
\tikzset{treeemph/.style={}}
\newcommand*\ATb{\begin{tikzpicture}[setree]\placeroots{1}\end{tikzpicture}}
\newcommand*\ATbb{\begin{tikzpicture}[setree]\placeroots{1}\children{child{node{}}}\end{tikzpicture}}
\newcommand*\ATbapab[1][]{\begin{tikzpicture}[setree, #1]\placeroots{2}\jointrees{1}{1}\end{tikzpicture}}

\def\d{\mathrm{d}}

\begin{document}
\title{Butcher series}\thanks{R.I.~McLachlan is supported by the Marsden Fund of the Royal Society of New Zealand.
K. Modin is supported by the Swedish Foundation for Strategic Research (ICA12-0052) and EU Horizon 2020 Marie Sklodowska-Curie Individual Fellowship (661482).
% General acknowledgments should be placed at the end of the article.
}
\subtitle{A story of rooted trees and numerical methods for evolution equations}
\titlerunning{Butcher series}        % if too long for running head
\author{Robert I.\ McLachlan \and
        Klas Modin \and
        Hans Munthe-Kaas \and
        Olivier Verdier
}
\authorrunning{McLachlan, Modin, Munthe-Kaas, and Verdier} % if too long for running head

\institute{%
	R.I.\ McLachlan \at
	Institute of Fundamental Sciences, Massey University, New Zealand \\
	% Tel.: +123-45-678910\\
	% Fax: +123-45-678910\\
	\email{r.mclachlan@massey.ac.nz}
	\and
	K.\ Modin \at
	Mathematical Sciences, Chalmers University of Technology and University of Gothenburg, Sweden \\
	% Tel.: +123-45-678910\\
	% Fax: +123-45-678910\\
	\email{klas.modin@chalmers.se}
	\and
	H.\ Munthe-Kaas \at
	Department of Mathematics, University of Bergen, Norway \\
	% Tel.: +123-45-678910\\
	% Fax: +123-45-678910\\
	\email{hans.munthe-kaas@uib.no}
	\and
	O.\ Verdier \at
	Department of Computing, Mathematics and Physics, Western Norway University of Applied Sciences \\
	% Tel.: +123-45-678910\\
	% Fax: +123-45-678910\\
	\email{olivier.verdier@hvl.no}
}

\date{\today \\ To appear in \emph{Asia Pacific Mathematics Newsletter}}

\maketitle

\begin{abstract}
	%%
	%% One or two sentences providing a basic introduction to the field.
	%%
	Butcher series appear when Runge--Kutta methods for ordinary differential equations are expanded in power series of the step size parameter.
	%They provide an indispensable tool for determining the order of convergence of a method.
	%%
	%% Two to three sentences of more detailed background, comprehensible to scientists in related disciplines.
	%%
	Each term in a Butcher series consists of a weighted elementary differential, and the set of all such differentials is isomorphic to the set of rooted trees, as noted by Cayley in the mid 19th century. %; he found trees useful for enumerating alkanes.
	A century later Butcher discovered that rooted trees  can also be used to obtain the order conditions of Runge--Kutta methods, and he found a natural group structure, today known as the Butcher group.
	% This group structure was rediscovered by Connes and Kreimer in the late 1990s in the context of renormalization in quantum field theory.
	It is now known that many  numerical methods also can be expanded in Butcher series; these are called B-series methods.
	%%
	%% One sentence clearly stating the general problem being addressed in this particular study.
	%%
	A long-standing problem has been to characterize, in terms of qualitative features, all B-series methods.
	%%
	%% One sentence summarizing the main results (with the words "here we show" or their equivalent).
	%%
	Here we tell the story of Butcher series, stretching from the early work of Cayley, to modern developments and connections to abstract algebra, and finally to the resolution of the characterization problem.
	This resolution introduces  geometric tools and perspectives to an area traditionally
	explored using analysis and combinatorics.
	
	%%
	%% Two or three sentences explaining what the main results reveals in direct comparison to what was thought to be the case previously, or how the main result adds to previous knowledge.
	%%
%	Every B-series method is affine-equivariant, i.e., independent of the choice of affine coordinates, so a suggestion was that affine-equivariance might be a defining property.
%	The recent advances reported on here refute this suggestion and unlocks the solution: B-series methods are those that are equivariant with respect to affine transformations between vector spaces of all possible (finite) dimensions.
%	%%
%	%% One or two sentences to put the results into a more general context.
%	%%
%	This characterization of B-series methods opens up an entirely new perspective on a long-studied topic: traditionally it has been explored using analysis and combinatorics; now, geometry has entered.
	%%
	%% Two or three sentences to provide a broader perspective, readily comprehensible to a scientist in any discipline if that considerably enhances the paper.
	%%
	% This geometric insight 

	% \textbf{Keywords:} Butcher series, order conditions, numerical integrators, ordinary differential equations, rooted trees, elementary differentials, affine equivariance

	% \textbf{MSC2010:} 65-03, 01-08, 65L06

\end{abstract}

\section{From Cayley to Butcher}
Butcher series are mathematical objects that were introduced by 
the New Zealand mathematician John Butcher in the 1960s. He
introduced them as part of his study of Runge--Kutta methods,
a popular class of numerical methods for evolution equations such
as initial-value problems for ordinary differential equations, and
they remain indispensable in the numerical analysis of differential
equations.  In this article we provide
a brief introduction to Butcher series, survey their early history
up to their introduction by John Butcher, and relate the
story of the many connections that have recently been discovered
between Butcher series and other parts of mathematics,
notably algebra and geometry.%
\footnote{This article is not a comprehensive review and is focussed
on our own interests. Useful companions to this article are the detailed mathematical review
of Butcher series by Sanz-Serna and Murua \cite{sa-mu} and the textbook treatments
of Hairer et al. \cite{hlw,ha-no-wa}.}
We begin, however, with the traditional definition.

Butcher series are intimately associated with the set of smooth (infinitely differentiable) vector fields on vector spaces.
Indeed, let $f$ be a smooth vector field on a vector space $V$, % equipped with basis coordinates $x=(x^1,\ldots,x^n)$.
defining the ordinary differential equation (ODE)
\begin{equation}
\label{eq:ode} \dot x = f(x),
\end{equation}
where $\dot x =  \frac{\d x}{\d t}$ denotes the derivative with respect to time~$t$.
One way to study (\ref{eq:ode}) is to  develop the Taylor series of its solutions. 
Let $x(h)$ be the solution to (\ref{eq:ode}) at time $t=h$ subject to the initial condition $x(0)=x_0$.
The Taylor series of $x(h)$ in $h$ is 
\begin{equation}\label{eq:taylor_series}
x(h) = x(0) + h \dot x(0) + \frac{1}{2}h^2 \ddot x(0) + \dots.	
\end{equation}
We already know that $x(0)=x_0$ and $\dot x(0) = f(x_0)$. %, so that $(\dot x)|_{h=0}=f(x_0)$.
The additional terms can be found by repeatedly applying the chain and product rules.
For example,
$$\ddot x = \frac{\d}{\d t}\dot x = \frac{\d}{\d t} f(x) = f'(x)\dot x = f'(x) f(x),$$
or, relative to a basis in which $x=x^1\mathbf{e}_1 + \ldots+ x^n\mathbf{e}_n$, %$x^i$ on $V$, 
$$ \ddot x^i = \sum_{j=1}^{n} \frac{\partial f^i}{\partial x^j}(x) f^j(x),$$
where $f(x) = f^1(x)\mathbf{e}_1 + \ldots + f^n(x)\mathbf{e}_n$.
Continuing in this way gives
\begin{equation}\label{eq:eldiffexp}
\begin{split}
\dot x &= f(x),\\
\ddot x &= f'(x)f(x),\\
\dddot x &= f'(x)f'(x)f(x) + f''(x)(f(x), f(x)),\\
\ddddot x &= f'(x)f'(x)f'(x)f(x) + f'(x)f''(x)(f(x),f(x)) + \\
& \qquad 3 f''(x)(f'(x)f(x),f(x)) + f'''(x)(f(x),f(x),f(x)),\\
&\;\,\vdots 
\end{split}
\end{equation}
Here the $k$th derivative $f^{(k)}(x)$ of the vector field $f$ is regarded as a multilinear map $V^k\to V$.
For example, $f''(f,f)$ is the vector field on $V$ whose $i$th coordinate is
$$ \sum_{j,k = 1}^n \frac{\partial^2 f^i}{\partial x^j \partial x^k}(x)f^j(x) f^k(x).$$
A vector field of the form appearing in~\eqref{eq:eldiffexp}, combining $f$ and its derivatives, is called an \emph{elementary differential}.
Using~\eqref{eq:eldiffexp}, the Taylor series~\eqref{eq:taylor_series} for the solution of (\ref{eq:ode}) can be written as 
\begin{equation}\label{eq:taylor_exact}
x(h) = x_0 + h f + \frac{1}{2}h^2 f'f + \frac{1}{6}h^3 f'f'f + \frac{1}{6}h^3 f''(f,f) + \dots
\end{equation}
where each elementary differential is evaluated at~$x_0$.
Notice that the power of $h$ in each term is determined by the multiplicity of $f$ in the elementary differential.
% : the power of $h$ is determined by the elementary differential.
However, the coefficients $1$, $1$, $1/2$, $1/6$, $1/6$, and so on are \emph{not} determined by their corresponding elementary differentials.
A \emph{Butcher series}, shortly denoted \emph{B-series}, is a generalization of~\eqref{eq:taylor_exact} allowing arbitrary coefficients, i.e., a formal series of the form
% Thus, formal series
% the Butcher series with coefficients $c = (c_0,c_1,\ldots)$ is given by 
% it is often convenient to take $h=1$ and instead take the step size as the scaling of~$f$.
% % That is, as a weighted sum of elementary differentials.
% Butcher series are generalizations to 
% This is an instance of the  
% series of the form
\begin{equation}
\label{eq:bseries}
 B(c,f) := c_0 x_0 + c_1 h f + c_2 h^2 f'(f) + c_3 h^3 f'(f'(f)) + c_4 h^3 f''(f,f) + \dots
 \end{equation}
where~$c_i\in\RR$. %are real coefficients. 
% The new vector fields $f'(f)$ and so on, called {\em elementary differentials} of $f$, are all evaluated at the same base point $x_0$. 
% A Butcher series, or B-series, is a formal series of the form (\ref{eq:bseries}).
Although presented here in coordinates, we shall see that Butcher series do not depend on the choice of basis.

\section{Early history}

Butcher series are named in honour of the New Zealand mathematician John Butcher. In a publication career spanning (so far) 60 years he has written 167 papers and books, all but 18 of them concerned  with Runge--Kutta methods and their generalisations. Most of them involve in some way  the fundamental structure that bears his name. Butcher series were introduced in a remarkable series of ten sole-authored papers in the years 1963--1972. 

A Runge--Kutta method is a numerical approximation $x_{n} \mapsto x_{n+1}$ of the exact flow of~\eqref{eq:ode}
% such that $x_k$ approximates the solution of the ordinary differential equation $\dot x = f(x)$ at $t$
defined by the following equations 
in $x_n$, $x_{n+1}$, $X_1,\dots,X_\nu\in V$:
\begin{equation}
\label{eq:rk}
\begin{aligned}
X_i &= x_n + h \sum_{j=1}^\nu a_{ij} f(X_j), \\
x_{n+1} &= x_n + h \sum_{j=1}^\nu b_j f(X_j).\\
\end{aligned}
\end{equation}
Here $\nu$ is the number of stages of the method and $a_{ij}$, $b_j$ are real numbers parameterising the Runge--Kutta method. Associated with the abstract Runge--Kutta method (\ref{eq:rk}) are its {\em order conditions}, polynomials equations in $a_{ij}$ and $b_j$---one equation per
elementary differential---that determine the order of convergence of the method and
its local error. Their derivation has been simplified over the years; a modern exposition
can be found in Hairer, Lubich and Wanner \cite{hlw}, and a detailed history in Butcher and Wanner \cite{bu-wa}.

The first breakthrough paper dates from 1963 \cite{butcher63}. Here Butcher found for the first time the coefficients $c_i$ of the B-series (\ref{eq:bseries}) of $x_{n+1}$ of the  Taylor expansion  in $h$ of an arbitrary Runge--Kutta method. This gave the order conditions for Runge--Kutta methods in complete generality. As  previous studies had laboriously expanded the solutions of particular (e.g. explicit) methods by
hand, this was an enormously important development. 

Butcher did have, however, some precursors. The most notable example is the paper of Merson~\cite{merson} from 1957. Robert Henry `Robin' Merson (1921--1992) was a scientist at the Royal Aircraft Establishment, Farnborough, UK, who was invited along with  more senior numerical analysts to a conference on Data Processing and Automatic Computing Machines at Australia's Weapons Research Establishment in Salisbury, South Australia.%
\footnote{%
Flight-related research at Farnborough began with the Army Balloon Factory in 1904, which
became the Royal Aircraft Factory in 1912, the Royal Aircraft Establishment in 1918,
and then the Royal Aerospace Establishment in 1988. It was merged into the 
Defence Research Agency in 1991 and then into the Defence Evaluation and Research
Agency in 1995. This was split up in 2001, with  Farnborough  becoming
part of the private company Qinetiq. Desmond King-Hele's version of these later 
developments is recorded at \cite{kh}.}
 It seems like a long way to go for a conference in 1957. However, the UK was still performing above-ground atomic bomb tests in South Australia at that time and the Australian government was very keen to be a part of the emerging era. Merson's work is bound up with one of the most significant events of 1957, the launch of Sputnik 1 on 4 October 1957, and the tale of Farnborough's involvement is  told in detail by one of the key participants, Desmond King-Hele, in his book {\em A Tapestry of Orbits} \cite{heleking}. The short version is that with the aid of a large radio antenna hastily erected in a nearby field, and some calculations of Robin Merson, within two weeks they had an accurate orbit for Sputnik 1. This allowed them to
estimate the density of the upper atmosphere and (after Sputnik 2) the shape of the earth.
Robin Merson became an expert in practical numerical analysis and orbit determination.

Merson's paper explains clearly the structure of the elementary differentials $f'(f)$, $f''(f,f)$, etcetera, and, crucially, shows how they are in one-to-one correspondence with rooted trees. He also introduces  various basic operations on rooted trees. This development, perhaps regarded initially as a bookkeeping device for finding and keeping 
track of the different terms, has over time become central to the combinatorial and algebraic study of B-series.

\begin{figure}
\begin{center}
\includegraphics[width=8cm]{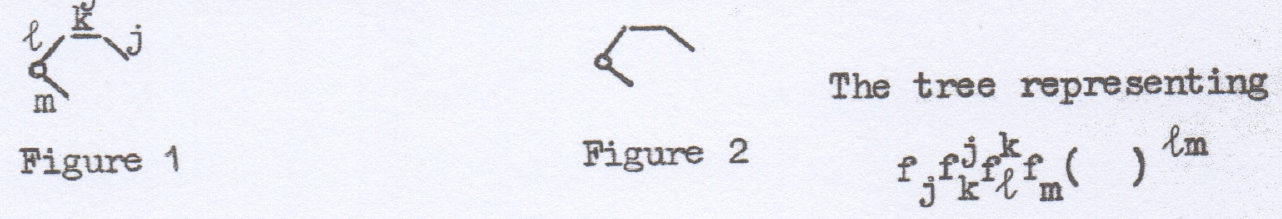}\\[4mm]
\includegraphics[width=8cm]{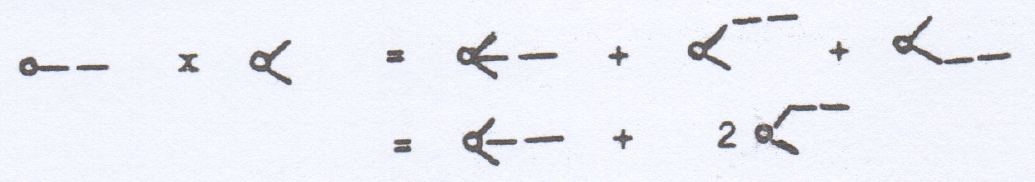}
\caption{\label{fig:merson}
Merson's \cite{merson} 1957 diagram of rooted trees representing elementary differentials, and (bottom) an example
of a product of trees, in this case the pre-Lie product explained in Section \ref{sec:alg}.
}
\end{center}
\end{figure}

\begin{figure}
\begin{center}
\includegraphics[width=8cm]{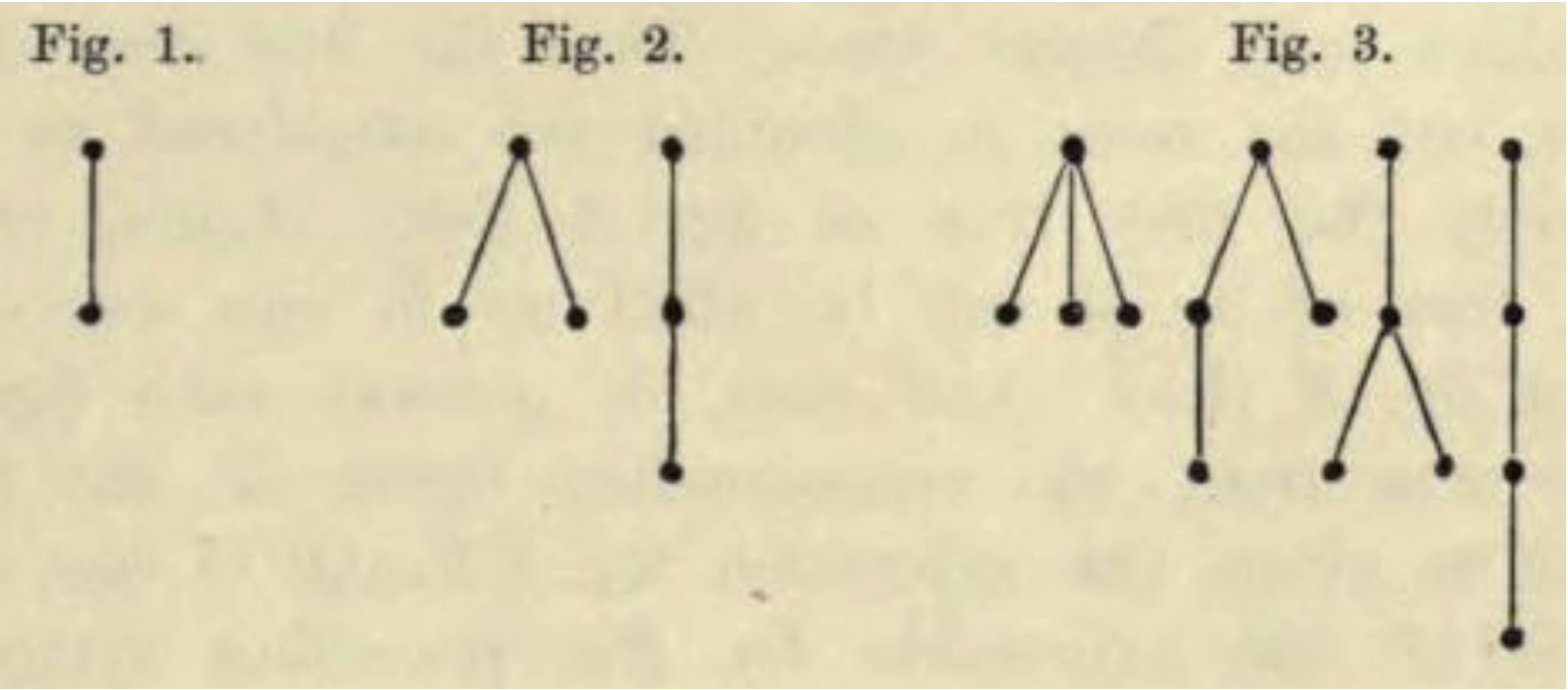}
\caption{\label{fig:cayley}
Cayley's \cite{cayley} 1857 diagram of rooted trees representing elementary differentials.
}
\end{center}
\end{figure}

The rooted trees $\mathcal{T}$ and their associated elementary differentials $\mathcal{F}(\mathcal{T})$ are 

\setlength{\tabcolsep}{0pt}
\noindent
\begin{tabular}{rlclclclclclclclclc}
$\mathcal{T} =\Big\{\emptyset$&,
&$\ab$&,
&$\aabb$&,
&$\aaabbb$&,
&$\aababb$&,
&$\aabababb$&,
&$\aabaabbb$&,
&$\aaababbb$&,
&$\aaaabbbb$&,
&$\dots
\Big\}$,
\\[4mm]
$\mathcal{F}(\mathcal{T})  = 
\Big\{
x$&,\; 
&$ f$&,\;
&$ f'(f)$&,\;
&$f'(f'(f))$&,\;
&$f''(f,f)$&,\;
&$f'''(f,f,f)$&,\;
&$f''(f,f'(f))$&,\;
&$f'(f''(f,f))$&,\;
&$f'(f'(f'(f)))$&,\;
&$\dots\Big\}.$
\end{tabular}
\setlength{\tabcolsep}{2pt}

\medskip

% \begin{align*}
% 	\mathcal{T} &= \Big\{\emptyset, &\ab \\
% 	\mathcal{F}(\mathcal{T}) &= \Big\{x, & f 
% \end{align*}

Merson introduces a method
for carrying out the required Taylor series expansions in elementary differentials and gives an example
of a 4th order Runge--Kutta method he derived. However, the actual expansions, although greatly simplified by the use of elementary differentials and rooted trees, are still carried out term by term. He did not have the coefficients of
all elementary differentials at once, as Butcher achieved. 

As it happens, the required mathematics and structures had already been discovered a century earlier by Arthur Cayley in 1857 \cite{cayley} (see Fig.~\ref{fig:cayley}). This is the actual discovery of the objects called trees (connected, cycle-free graphs). In popular treatments of graph theory, the development of graph theory is closely linked with recreational mathematics (the bridges of K\"onigsberg) and with chemistry (Cayley's enumeration of alkanes and other families of molecules). One common interpretation of the story is that Cayley introduced the trees as a purely abstract structure and 17 years later---behold the power of mathematics!---found that he could use them to count molecules. However, Cayley actually needed trees for {\em exactly} the purpose we are using them here---to keep track of how vector fields interact when applied repeatedly to one another---and this purpose was then forgotten for a hundred years. 
As the need for better numerical integration methods arose towards the end of the 19th century, the required tools for a complete theory were indeed already there, but they had been forgotten.

As Frank Harary wrote \cite{harary}, 

{\narrower\narrower\medskip

\noindent\em 
In very many cases and in disciplines in the physical sciences, the social sciences, computer science, and the humanities, graphs frequently occur as a natural, useful, and intuitive mathematical model. The consequence is
that those investigators who were not aware of the existence of graph theory as a study in its own right were led to rediscover it in order to apply it.

\medskip
}

\noindent
Interestingly enough, Merson does cite Cayley. However, from the context, it is not clear that he actually laid eyes on Cayley's paper. 
He writes,

{\narrower\narrower\medskip

\noindent\em A formula for the number of trees of a given order was discovered by CAYLEY {\rm [our \cite{cayley}]}
and quoted by ROUSE--BALL\dots

\medskip
}

\noindent
This was probably the original 1892 edition of Rouse Ball's famous book {\em Mathematical Recreations and Essays}, as later
editions included Coxeter as coauthor. This first edition contains just one page on trees, stating
Cayley's formulae for the number of trees. Now this same section of Rouse Ball also discusses the famous Knight's Tour problem, an astonishingly long-lived problem dating from an Arabic manuscript of 840 AD. For example, there were three articles on Knight's Tours published in the {\em Mathematical Gazette} in 1956 alone. This problem became a life-long interest of Merson's, who published tours in 1974 and 1999 (posthumously, in {\em Games and Puzzles} magazine, from letters written in 1990--91)
that are still in many cases the best known tours. Although Merson
stated \cite{jelliss} that he first became interested in the 
problem in 1972,  it is not unlikely that in 1957 he rediscovered trees independently because, like Cayley,  he needed them, and from his interest in recreational mathematics remembered Rouse Ball's discussion of Cayley without ever chasing it up.

John Butcher, at that time a PhD student in physics at the University of Sydney, was actually present at Robin Merson's talk in 1957, but says \cite{jcbpc} that
he did not understand it at all. However,  the seed was planted there.
To return to Butcher's 1963 paper, he closes with the following statement: 

{\narrower\narrower\medskip

\noindent\em It happens that this situation is capable of extensive generalization and, for example, keeping this same value $\nu=3$ it is possible to satisfy the 37 conditions necessary for a sixth order process. Similarly for any value of $\nu$ a process of order up to $2\nu$ is possible. It is intended that details of such processes will be discussed in a later publication.

\medskip
}

This was an announcement of Butcher's discovery of the family of Gauss Runge--Kutta methods and
the first hint of extra structure contained within the Runge--Kutta order conditions.
Methods with 3 stages have 12 free parameters ($a_{ij}$ and $b_j$ for $i,j=1,2,3$) and Butcher was extremely excited to discover that
there were values of the parameters that satisfied not just the 8 conditions for order 4,
and the 17 conditions required for order 5, but even the 37 conditions required for order
6! He recalls running through the empty corridors of the mathematics department at the University
of Canterbury, where he was then lecturing, desperately trying to find someone to
understand and to share the excitement \cite{jcbpc}. He fulfilled his intention to publish the details
in his very next paper \cite{butcher64}.

One approach taken by Butcher to approach the structure of the order conditions, suggested by this discovery, was to introduce certain {\em simplifying assumptions}. These became the cornerstone of the construction of the efficient high-order explicit integrators that are used today. However, the source of these
simplifying assumptions  remained mysterious; only very recently has their algebraic origin been 
explained \cite{khashin}. This has allowed them to be embedded in systematic families and further reduced
the number of stages needed at high order. We take this as further evidence that after 50 years Butcher's vision is alive and well.

This initial intensely creative and productive period came to a head with the publication of {\em An algebraic theory of integration methods} in 1972 \cite{butcher72}---submitted in 1968---in which John Butcher introduced what is
now called the Butcher group. The B-series (\ref{eq:bseries}) with $c_0=1$ correspond formally to
diffeomorphisms close to the flow of $f$, and the Butcher group operation arises from a product
of rooted trees that corresponds to the composition of these diffeomorphisms.

To give an example of the group operation of the Butcher group, consider the B-series
$$\alpha := x_0 + h f(x_0).$$
This is associated with the map $x_0 \mapsto x_1 := x_0 + hf(x_0)$ of the forward Euler method.
The composition of this map with itself (i.e., two steps of forward Euler)
 is the map 
 \[
 \begin{aligned}
 x_0& \mapsto x_1 + h f(x_1) \\
 &= x_0 + h f(x_0) + h f(x_0 + h f(x_0))\\
 &= x_0 + h f + h(f + h f'f + \frac{1}{2!} h^2 f''(f,f) + \frac{1}{3!}h^3 f'''(f,f,f) + \dots)\\
 &= x_0 + 2 h f + h^2 f'f + \frac{1}{2!}h^3 f''(f,f) + \frac{1}{3!} h^4 f'''(f,f,f) + \dots.
 \end{aligned}
 \]
 The last line is the B-series of the Butcher product $\alpha\alpha$.
 
 The inverse $\alpha^{-1}$ of the B-series $\alpha$ is the series associated with the inverse map $x_1\mapsto x_0$.
 This map is one step of backward Euler with time step $-h$. Its B-series is
 \[
 \begin{aligned}
 x_0 - h f & + h^2 f'f - h^3(f'f'f + \frac{1}{2}f''(f,f)) \\
 & + h^4(\frac{1}{6}f'''(f,f,f) + f'f'f'f + f''(f,f'f) + \frac{1}{2}f'(f''(f,f))) + 
\dots.
\end{aligned}
\]
The coefficient of any elementary differential in these series can be found using simple
combinatorial operations on trees.

This paper \cite{butcher72} aroused an interest that lead to a crucial event.
In Innsbruck, the 28-year-old dozent Gerhard Wanner was studying John Butcher's early papers and 
his hard-to-understand preprint \cite{butcher72}. In 1970 the University of Innsbruck was celebrating
its 300th anniversary and asked each professor to invite a guest lecturer. Wanner's professor,
Wolfgang Gr\"obner, asked Wanner for a suggestion, and so John Butcher was invited. Ernst Hairer,
who had been Wanner's best freshman analysis student the year before, attended the lectures. In Wanner's words \cite{wanner}, \emph{%
``In my opinion, at that time, nobody in the world made the necessary efforts to understand Butcher's papers, except Ernst. He then explained them to me, and I tried to put them in a more understandable form,"} and in 
Butcher's words~\cite{butcherearly}, \emph{``This led to my own contribution being recognised, through their eyes, in a way that might otherwise not have been possible."}
In 1974 Hairer and Wanner \cite{ha-wa} introduced both
 {\em Butcher series} and the term {\em Butcher group}; they also clearly
demonstrate the uses of the series for much more than Runge--Kutta methods.
In Butcher \cite{butcher72}, the group elements are functions from rooted trees to the reals,
such as those functions induced from (traditional and continuous stage) Runge--Kutta methods;
in Hairer and Wanner \cite{ha-wa} the
primary objects are  the B-series (\ref{eq:bseries}) themselves, which obey the  group law  found by Butcher.

These discoveries triggered a period of huge development in numerical methods
for evolution equations. 
The subsequent modern history of the area has been reviewed extensively \cite{bu-wa,hlw,ha-no-wa,sa-mu}.
Here we confine ourselves to some remarks as to the role and significance of Butcher series.

\section{How important are Butcher series?}

Many areas of inquiry show a tendency to divide adherents into `lumpers' and `splitters'. For 
example, in taxonomy, lumpers prefer to name few species, splitters many. Lumpers 
emphasize similarity, splitters emphasize difference. Numerical analysis,
like most parts of mathematics, shows a gradual tendency over time towards splitting, as 
the true differences between instances are appreciated and exploited. 
Thus structure-preserving
methods have been developed for finer and finer divisions of matrices, differential equations and
so on, that, by restricting the problem class, are able to offer superior performance. 
Iserles \cite{iserles} alludes to this when he compares ordinary differential equations to Tolstoy's happy families, that (`perhaps', Iserles cautions) all resemble each other, while each partial differential 
equation is unhappy in its own way. Indeed, a mighty strength, and also a potential
weakness, of Runge--Kutta methods and of B-series is that they treat {\em all} ODEs in a uniform way. They are an extreme example of lumping. One might wonder if they are perhaps {\em too} extreme.
Do they over-lump ODEs?

In our view they have held up pretty well. The first widely-acknowledged division of ODEs in numerical analysis was into stiff and nonstiff equations. Implicit Runge--Kutta methods turned out to be ideal for stiff equations and explicit ones for nonstiff. With the advent of symplectic integrators for Hamiltonian systems, that preserve a quadratic conservation law on first variations of solutions, Runge--Kutta methods were found to be suitable too. New classes of methods have been introduced that have features that Runge--Kutta methods do not, such as exponential integrators like
\begin{equation}
\label{eq:ei}
x_{n+1} = x_n + \phi(h f'(x_n))h f(x_n),\quad \phi(z) = \frac{e^z-1}{z},
\end{equation}
which can beat implicit Runge--Kutta methods on some stiff equations, and the AVF (Average Vector Field)
method
\begin{equation}
\label{eq:avf}
x_{n+1} = x_n + \int_0^1 f(\xi x_{n+1} + (1-\xi) x_n)\, {\rm d}\xi
\end{equation}
that preserves energy $H(x)$ when $f = J^{-1}\nabla H$ is a Hamiltonian vector field. 
Both (\ref{eq:ei}) and (\ref{eq:avf}) have expansions in B-series.

On the other hand, some methods such as the leapfrog or St\"ormer--Verlet
method, widely used in molecular dynamics and in video game engines for systems of the form $\ddot x = -\nabla V(x)$, do not have B-series---indeed they are not even defined for all first
order systems $\dot x = f(x)$---and should certainly not be discarded on that account. Our view
is {\em lump if you can, but split if you must}. 

In fact some
would say that there is {\em no} practical reason for preferring methods with a B-series and that
the whole concept is merely a mathematical abstraction or (perhaps) convenience. 
However, note that (\ref{eq:bseries}) lumps not only ODEs, but also numerical methods. A very large
class of numerical methods for ODEs are represented by (\ref{eq:bseries}). 
Even before getting to the question of what the possession of a B-series confers on a numerical method, the lumping of numerical methods by B-series presents a fairly rare opportunity in 
computational science. All too often one analyzes the complexity or behaviour of a {\em particular} algorithm, or perhaps of a small class. Meaningful lower bounds for complexity or behaviour over {\em all} algorithms are almost never obtained. One should not miss the opportunity given by B-series to better understand an infinite-dimensional set of methods, without regard to particular details of the method.

Several times, new numerical methods have been reflected in the discovery of new structure
within B-series. For example, if $f=J^{-1}\nabla H$ for some $H$ and $J$, where $J^T=-J$
defines a symplectic structure on the vector space $V$, then $f$ is Hamiltonian and
energy preserving and we can ask which B-series have these properties. 
The trivial B-series $B(f)=c_1 f$ 
are the only ones which are both Hamiltonian and energy-preserving.
At first sight it
is surprising that the first nontrivial B-series, $f'f$, is neither Hamiltonian nor energy-preserving. At the next order, $f'f'f$ is energy preserving and $f''(f,f)-2f'f'f$ is Hamiltonian. The spaces of such
B-series have been completely described \cite{cmoq}.

\section{Algebraic characterizations}
\label{sec:alg}

The topic of B-series can be approached from many different points of view; topics in numerical analysis, geometry and abstract algebra are connected via B-series. The fundamental algebraic structure of a \emph{pre-Lie algebra} 
unifies three seemingly very different papers all written in 1963: John Butcher's first paper on Runge--Kutta methods~\cite{butcher63}, Ernest Vinberg's paper on the geometry of symmetric cones~\cite{vinberg}
and Murray Gerstenhaber's work on homology and deformations of algebras~\cite{gerstenhaber63}. The differential geometric picture starts with the basic notion of parallel transport of vectors, which is infinitesimally described in terms of a \emph{connection} or \emph{covariant derivation} of vector fields.
The connection is a bilinear operation of vector fields $(f,g)\mapsto f\tr g$ (often written as $\nabla_f g$) which  describes the 
rate of change of~$g$ as it is parallel-transported along the flow of~$f$. On the vector space $\RR^n$ parallel transport is the obvious rule, and the corresponding connection 
is given as 
\[f\tr g = g'(f) = \sum_{i,j=1}^n \frac{\partial g^i}{\partial x^j} f^j\ \frac{\partial}{\partial x^i} .\]
The curvature  $R$ and the torsion  $T$ are the two basic invariants of a connection. On flat spaces, such as the above defined connection on $\RR^n$, both $R=0$ and $T=0$. It can be shown that in this case the connection satisfies the following \emph{pre-Lie} relation:
\[f\tr(g\tr h)- (f\tr g)\tr h =g\tr(f\tr h)- (g\tr f)\tr h.\]
An algebra with a product satisfying this relationship is called a \emph{pre-Lie algebra}. So, the set of smooth vector fields on $\RR^n$ with the standard connection is an example of a pre-Lie algebra\footnote{Also called a Vinberg, Koszul--Vinberg, left-symmetric, or Gerstenhaber algebra. The name reflects the fact that the skew product $[x,y] := x\tr y-y\tr x$ defines a Lie bracket. However it should be noted that the pre-Lie relation is not the most general form of a product with this property.}. 
Another example is the linear combination of rooted trees, where the pre-Lie product is given by grafting: for two trees $\tau_1$ and $\tau_2$ the pre-Lie product $\tau_1\tr\tau_2$ is computed by attaching the root of $\tau_1$ with an edge to each of the nodes of $\tau_2$ and adding all these terms together (see Figure \ref{fig:merson}.) 
The pre-Lie algebra perspective of B-series was promoted by Calaque, Ebrahimi-Fard, and Manchon~\cite{kurusch_et_al}.
A fundamental result, which was essentially known already to Cayley in 1857, but which has been revisited in a modern algebraic setting by Chapoton and Livernet in 2001 \cite{ch-li}, is that the space of all trees with the grafting product is the \emph{free pre-Lie algebra}.
This means that this structure `knows all there is to know' about basic algebraic properties of pre-Lie algebras, and any algebraic computation which relies only on the pre-Lie relationship
can be expressed as a computation on trees. It also means that any example of a concrete pre-Lie algebra can be realised as a quotient of the free pre-Lie algebra with some ideal (that is, as trees with some equivalence relation). 
This is indeed a useful result for computations.

The correspondence between  abstract trees and  concrete elements in a given pre-Lie algebra (e.g.,\ a vector field on $\RR^n$) is exactly the elementary differential map of Butcher. 
The elementary differential map $\F(\tau)$, taking trees to vector fields, respects the structure of the pre-Lie product, $\F(\tau_1\tr\tau_2) = \F(\tau_1)\tr\F(\tau_2)$, where the triangle on the left is grafting of trees and on the right is the covariant derivative of vector fields. All the elementary differentials are obtained this way. 
For example, since
$\aababb = \ab\tr (\ab\tr\ab) - (\ab\tr\ab)\tr\ab$, we must have that if $\F(\ab) = f$, then
$\F(\aababb) = f\tr(f\tr f)-(f\tr f)\tr f$. 
Similarly, all the terms of the B-series can be expressed in terms of the pre-Lie product, and hence we can regard a B-series as an infinite expansion in a pre-Lie product. 

Are there other important examples of pre-Lie algebras where B-series might play a role? There was a great surprise in the late 1990s when Christian Brouder pointed out \cite{brouder} that the so-called Hopf algebra of Alain Connes and Dirk Kreimer \cite{co-kr} had the same algebraic structure that John Butcher had been studying in detail in his 1972 paper. Connes and Kreimer had been interested in renormalisation processes in quantum field theory and discovered a rich algebraic structure of trees. Indeed Arne D\"ur \cite{dur} had already observed in 1986 that Butcher had  given rooted trees the structure of a Hopf algebra.  Rereading
Butcher \cite{butcher72} in light of these more recent developments, it is striking how close his 
perspective is to the modern Hopf algebraic view. As Brouder commented, \emph{``Butcher found an explicit
expression for all the operations of the Hopf structure of the algebra of rooted trees.''}
After Brouder's work the Fields medallist Alain Connes wrote  \cite{co-kr} \emph{``We regard Butcher's work on the classification of numerical integration methods as an impressive example that concrete problem-oriented work can lead to far-reaching conceptual results.''} Pierre Cartier has also 
written a very clear exposition of the significance of pre-Lie algebras and the 
algebraic origin of the Connes--Kreimer approach \cite{cartier} .

More recently these algebraic structures appear in other important areas, such as in stochastic processes, where the \emph{Rough Paths Theory} gives a precise meaning to integrating functions along highly irregular paths. This theory originated from the work of Terry Lyons and was celebrated by the Fields medal awarded to Martin Hairer in 2014 for his work on regularity structures.
Relations between rough paths and B-series have been developed in the work of Massimo Gubinelli \cite{gubinello}.

In a completely different direction, expansions in rooted trees can be used to dramatically simplify and also to sharpen known results in complex dynamics \cite{fauvet} ({\em ``this amounts to a novel approach to formal linearization by means of a powerful and elegant combinatorial machinery''}).

Considering B-series as an expansion in a (flat and torsion free) connection, we may ask what are the characterising geometric properties of a B-series? A partial answer comes from the
question of which invertible mappings $\phi\colon \RR^n\rightarrow \RR^n$ preserve the connection $\tr$. Let $\phi$ act on  vector fields in the `natural' way (i.e.,\ as a differential equation transforms under change of coordinates) $\phi\cdot f := (\phi')\circ f\circ \phi^{-1}$, where $\phi'$ is the Jacobian matrix. Then it can be shown that $\phi\cdot(f\tr g) = (\phi\cdot f)\tr (\phi\cdot g)$ for all vector fields $f$ and $g$ if and only if $\phi(x) = Ax+b$ is an affine map. However, it turns out that this condition is not enough to nail precisely the question of \emph{What is a B-series?}, but we shall see %in the next section
that it brings us a long way towards the answer.
Before we explore this issue further in the next section, we remark on other recent geometric developments of the theory. 

Concerning the group structure of B-series, Bogfjellmo and Schmeding~\cite{bogfjellmo_schmeding} have recently proved that the space of B-series is an infinite-dimensional Lie group with respect to a natural Fréchet topology.
Among numerical analysts, B-series have long been treated as Lie groups without a rigorous justification;
the result by Bogfjellmo and Schmeding resolves this and unveils interesting possibilities to apply tools from infinite-dimensional geometry to the backward error analysis of ODE methods.

The question of characterising geometries by invariance properties goes a long time back to the 19th century work of Felix Klein, who in his \emph{Erlangen program} of 1872 raised fundamental questions about geometries and symmetries. An example is the study of affine geometries as a generalisation of Euclidean spaces. In this geometric context it is interesting to ask if other geometries have algebras describing their connections, such as pre-Lie algebras for affine geometries. Recent developments have shown that this is indeed the case. For Lie groups and homogeneous spaces there are naturally defined connections which give rise to \emph{post-Lie} algebras, and from this we obtain B-series types of expansions valid for flows evolving on manifolds (`Lie--Butcher' series) \cite{mk-lu}. Yet another algebra appears in the context of  symmetric spaces such as, for example, spheres and Riemannian spaces with constant curvature.
This is an active area of research, where differential geometry, algebraic combinatorics, differential equations, computations and applications go hand-in-hand.

%Thus, the study of Butcher series is very definitely alive 50 years after their discovery. 
%Recently we have established a characterization of Butcher series that, at least to our minds, answers
%the title question, {\em What are Butcher series, really?} 

\section{Geometric characterizations}

Many mathematical objects can be defined in different ways: axiomatically, constructively, or by characterizing their relationship to another, known, object. 
The original, and still the traditional, approach to Butcher  series \cite{hlw} is constructive.
It is motivated
by the Taylor series of the exact solution. It starts by constructing the rooted trees, most easily done recursively using the operation of adding a root to a forest (set of rooted trees). Then the elementary differentials are defined and associated to the rooted trees, and finally it is shown that
various objects (Runge--Kutta and other integration methods) can be expanded in Butcher series.
The algebraic approach of the previous section is axiomatic. However,
if we recall the origin of Butcher series in numerical analysis, and note that not all numerical
integrators have a Butcher series, it is natural to ask 
why  these particular combinations, $f''(f,f)$ and so on, keep coming up. What is special about them?
What geometric property characterises those numerical integrators that have a Butcher series?

A crucial clue is provided in the definition of Runge--Kutta methods, (\ref{eq:rk}). 
Apart from evaluation of $f$, these involve
only  scalar multiplication and addition---the defining operations of the vector space~$V$. 
This suggests that Runge--Kutta methods are defined intrinsically on $V$ and do not
depend on the choice of basis. 
Indeed, as already mentioned previously in the context of pre-Lie algebras, slightly more is true: Runge--Kutta methods (and
B-series) are {\em affine-equivariant}. % in the following sense.
Indeed, let, as before, smooth invertible mappings $\phi\colon V\to V$ act
on the vector space $V$ and on vector fields on $V$ in the natural way. Then B-series with $c_0=1$,
such as the expansions of numerical integrators, obey
$$\phi\cdot B(c,f) = B(c,\phi\cdot f)$$
for  all invertible affine maps $\phi(x) = A x + b$, $A\in {\mathbb R}^{n\times n}$, $\det A\ne 0$. 
% This is the approach taken by Chapoton \cite{ch-li} who constructs B-series as the free object generated
% from the pre-Lie product $(f,g)\mapsto g'(f)$ of vector fields, which as he notes is intrinsically
% defined on an affine space. 
Could it be the case that any affine-equivariant method has a Butcher series?
In other words, does affine-equivariance characterize B-series methods?

In \cite{mko}, two of us showed that this is not the case. There are many methods that 
are affine-equivariant but do not have Butcher series. 
The simplest example is the first-order method
\[
x_1 = x_0 + h f(x_0)(1 + h (\nabla\cdot f)(x_0)).
\]
Under an affine transformation $x \mapsto \phi(x)=Ax+b$, $f$ transforms to $A f\circ\phi^{-1}$,
and the Jacobian $f'$ transforms to $A (f'\circ\phi^{-1}) A^{-1}$. The divergence of $f$, namely
${\rm tr} f'$, transforms to $({\rm tr} f')\circ\phi^{-1}$, and the new term $f\, \nabla\cdot f$
transforms to $A (f\,\nabla\cdot f)\circ\phi^{-1}$---that is, it is affine equivariant.

It turns out that any affine-equivariant method 
can be expanded in terms of more general objects, the {\em aromatic series}. 
Combinatorically, these are represented by `aromatic trees', forests consisting of one rooted tree and any number
of directed graphs with one cycle (self-loops allowed). The name is suggested by
aromatic compounds, such as benzene, that contain cycles of  atoms. An aromatic series begins
\begin{equation*}
\label{eq:aroma}
\begin{aligned}
c_0 x & + c_1 h f \\
&+ h^2(c_2 f'f + c_3f \nabla \cdot f ) \\
&+ h^3(c_4 f''(f,f) + c_5 f'f'f + c_6 f(f\cdot\nabla(\nabla\cdot f)) + c_7
f'f\,\nabla\cdot f \\
& \qquad+ c_8 f (\nabla\cdot f)^2 + c_9 f {\rm tr}(f'^2)) \\
& + \dots
\end{aligned}
\end{equation*}
which may be represented as an element in the span of the aromatic trees
\begin{equation*}
\begin{aligned}
	&\ATb,\\
	& \ATbb, \quad \ATbapab[treeemph],\\
	&
	\aababb
	,\quad
	\begin{tikzpicture}[setree]
		\placeroots{1}
		\children[1]{child{node{} child{node{}}}}
	\end{tikzpicture}
	,\quad
\begin{tikzpicture}[setree, treeemph]
	\placeroots{2}
	\children[1]{child{node{}}}
	\jointrees{1}{1}
\end{tikzpicture}
,
\quad
\begin{tikzpicture}[setree]
	\placeroots{2}[.3]
	\children[1]{child{node{}}}
	\joinlast{2}{2}
\end{tikzpicture}
,
\quad
\begin{tikzpicture}[setree]
	\placeroots{3}
	\joinlast{2}{2}
	\joinlast{3}{3}
\end{tikzpicture}
,
\quad
\begin{tikzpicture}[setree]
	\placeroots{3}[.3]
	\jointrees{2}{3}
\end{tikzpicture}
,
\\
&\ldots
\end{aligned}
\end{equation*}
 There are clearly many more aromatic than rooted trees. The aromatic trees  of
 order $n$
 are in 1--1 correspondence with functions from $\{2,\dots,n\}$ to
 $\{1,\dots,n\}$, `forgetting the labels', that is, modulo permutations of $\{2,\dots,n\}$. (Here the element
 1  identifies the root.)
 For example, the aromatic tree
 \begin{equation*}
 \begin{tikzpicture}
 \begin{scope}[etree, scale=1.5]
 \placeroots{2}[.3]
 \children[1]{child{node(a){}}}
 \children[2]{child{node(b){}}}
 \joinlast{2}{2}
 \end{scope}
 \node[left] at (tree1) {$1$};
 \node[above right] at (tree2) {$4$};
 \node[above left] at (a) {$2$};
 \node[above right] at (b) {$3$};
 \end{tikzpicture}
 \end{equation*}
is associated with the
 function $2\mapsto 1$, $3\mapsto 4$, $4\mapsto 4$ and with the (generalized) elementary differential
$$
\sum_{i_1,i_2,i_3,i_4=1}^n
  f^{i_1}_{i_2} f^{i_2} f^{i_3} f^{i_4}_{i_3 i_4}\frac{\partial}{\partial x^{i_1}} = f'(f)\, (f\cdot \nabla(\nabla\cdot f)).$$
The numbers of such `shapes of partially defined functions' is given in sequence A126285 in the Online Encyclopedia
 of Integer Sequences and tabulated in Table \ref{tab:1}. The number of
 rooted trees, first evaluated by Cayley, are shown for comparison.
 The apparently terrifying numbers of rooted trees were tamed by Butcher.
 What will happen to the even more plentiful aromatic trees?

\begin{table}
	\begin{center}
		\setlength{\tabcolsep}{4pt}

		\begin{tabular}{cp{4ex}p{4ex}p{4ex}p{4ex}p{4ex}p{4ex}p{4ex}p{4ex}p{4ex}p{4ex}}
			\toprule
			$n$ & 1 & 2 & 3 & 4 & 5 & 6 & 7 & 8 & 9 & 10\\
			\midrule
			\# rooted trees & 1& 1& 2& 4& 9& 20& 48& 115& 286& 719 \\
			\# aromatic trees & 1&  2&  6& 16& 45& 121& 338& 929& 2598& 7261\\
			\bottomrule
		\end{tabular}
		\caption{ \label{tab:1} Enumeration of rooted and aromatic trees with up to 10 nodes.}
	\end{center}
\end{table}
 
The existence of the aromatic series shows that  affine-equivariance of a method is not enough to ensure that it can be expanded in 
a B-series. What else is needed? The second big clue is that
Runge--Kutta methods are defined without reference to the
dimension of the underlying vector space. It does not seem to play any
role at all. Clearly, at a minimum, the expansion of the
method in each dimension must have the same coefficients. But what
rules out the aromatic terms like $f\,\nabla\cdot f$?

The answer is that these terms do not respect {\em affine-relatedness}. Consider two
vector spaces $V$ and $W$ of possibly different dimension, together with
an affine map $\phi\colon V\to W$, $x\mapsto A x + b$. The vector fields
$f$ on $V$ and $g$ on $W$ are said to be $\phi$-related if $g(A x + b) = A f(x)$ 
 for all $x\in V$. 
B-series preserve affine-relatedness in the sense that for any affine $\phi$,
if $f$ and $g$ are $\phi$-related then $B(c,f)$ is $\phi$-related
to $B(c,g)$. 
In \cite{mmmv} we prove that this property characterizes B-series: {\em a numerical method has a Butcher series if and only if it preserves affine-relatedness.}

Preserving affine-relatedness has a fairly direct physical interpretation.
It means that the method is immune to changes of scale, such as changes of units.
It means that the method preserves invariant affine subspaces {\em automatically}, whenever
the system has any such.
It means that the method preserves affine symmetries, again automatically; the method
does not even have to `know' (or be told) that the system has the symmetries.
It means that the method leaves decoupled systems decoupled, again automatically.
All these properties are desirable when designing general-purpose ODE software.
Furthermore, we now see that many of the more subtle properties of B-series, originally discovered
through combinatorial analysis of trees, must in fact be a direct consequence of
affine-relatedness. 
Examples include special properties with respect to symplecticity, preservation of quadratic invariants, and
preservation of energy \cite{chartier} and non-preservation of volume \cite{is-qu-ts}.

The proof of the theorem on affine equivariance \cite{mko} relies  on some classical results in 
functional analysis and invariant theory. First it is established that the Taylor series in $f$
of an arbitrary map depends only on the derivatives of $f$, and that the terms of order $n$
are in fact a polynomial of degree $n$ in $f$ and its partial derivatives. Second, the
invariant polynomials that are functions of $f$ and its partial derivatives, whose values at $x_0$ are regarded now
as arbitrary symmetric tensors, are sought using the `invariant tensor theorem'. The conclusion at 
2nd order is that only $f^i f^j_i$ and $f^i f^{j}_j$ are equivariant, these giving the two aromatic
trees of order 2. At 3rd order, to the tensor $f^i f^j f^k$ the partial derivatives $j$ and $k$ can
be attached to any two of the factors, leading to the 6 aromatic trees of order 3.

The proof of the theorem on affine relatedness, characterizing B-series \cite{mmmv}, begins
with an arbitrary affine-related method. Since, in particular, it is affine-equivariant, it has an aromatic series.
Each aromatic tree containing loops is to be knocked out. For each such tree, a special pair
of affine-related vector fields is constructed such that affine-relatedness of the method means
that the coefficient of this tree must be zero. For example, for the tree $\ATbapab$, associated
with $f\,\nabla\cdot f$, the vector fields are $f^{(1)}\colon\dot x^1 = 1$, $\dot x^2 = x^2$ and $f^{(2)}\colon\dot x^1=1$.
These vector fields are related by the affine map $(x^1,x^2)\mapsto x^1$. 
Since $f^{(1)}\,\nabla\cdot f^{(1)}=1$ and $f^{(2)}\,\nabla\cdot f^{(2)}=0$,
this term cannot appear in the expansion of a method that preserves affine-relatedness.

To summarize, Butcher series are objects intrinsically associated to 
the set of vector fields on affine spaces of {\em all} dimensions, and will
show up naturally in any analysis that respects the affine structure and does not depend on the dimension. 
This explains their ubiquity. 
It is fascinating that natural and practical demands of numerical methods for ODE---black-box solvers defined uniformly on all affine spaces---has led to the discovery of a fundamental invariant object.

On the other hand, where does this leave the aromatic series? We suggest that
they will show up naturally in problems posed in a specific dimension. Although traces and divergences
are common in physics, we have not seen aromatic series before. They arose
purely from a question in numerical analysis, but are fundamental in their own way.
Moreover, they can have properties that no B-series can have.
For example, many aromatic series, but no B-series, are divergence free.

\medskip
\noindent{\bf Acknowledgements}
We thank John Butcher, Ernst Hairer, and Gerhard Wanner for their comments.


\begin{thebibliography}{99}

\bibitem{bogfjellmo_schmeding} Bogfjellmo, G. and Schmeding, A., The Lie group structure of the Butcher group, {\em Found.\ Comp.\ Math.} (2015), DOI:10.1007/s10208-015-9285-5
\bibitem{brouder} Brouder, C., Runge--Kutta methods and renormalization, {\em Eur. Phys. J. C \bf 12}  (2000),521--534.
\bibitem{jcb} {\tt http://jcbutcher.com/publications}
\bibitem{jcbpc} Butcher, J. C., personal communication.
\bibitem{butcher63} Butcher, J. C., Coefficients for the study of Runge-Kutta integration processes, {\em J. Austral. Math. Soc. \bf 3} (1963), 185--201.
\bibitem{butcher64} Butcher, J. C., Implicit Runge-Kutta processes, {\em Math. Comp. \bf 18} (1964), 50--64.
\bibitem{butcher72} Butcher, J. C., An algebraic theory of integration methods, {\em Math. Comp. \bf 26} (1972), 79--106.
\bibitem{butcherearly} Butcher, J. C., Numerical methods for ordinary differential equations: early days, in
{\em The Birth of Numerical Analysis}, 
A. Bultheel and  R. Cools, eds., 
World Scientific, 2010, pp. 35--44.

\bibitem{bu-wa} Butcher, J. C., and Wanner, G., Runge-Kutta methods: some historical notes, {\em Appl. Numer. Math. \bf 22} (1996), 113--151.

\bibitem{kurusch_et_al} Calaque, D., Ebrahimi-Fard, K., and Manchon, D., Two interacting Hopf algebras of trees: A Hopf-algebraic approach to composition and substitution of B-series, {\em Adv.\ Appl.\ Math. \bf 47} (2011), 282--308.
\bibitem{cartier} Cartier, P., Vinberg algebras, Lie groups and combinatorics, in {\em  Clay Mathematics Proceedings. Quanta of Maths \bf 11} (2010), 107--126.
\bibitem{cayley} Cayley, A., On the theory of the analytical forms called trees, {\em Philos. Mag. \bf 13}(85) (1857), 172--176.
\bibitem{ch-li} Chapoton, F. and Livernet, M, Pre-Lie algebras and the rooted trees operad, {\em International Mathematics Research Notices \bf 8} (2001), 395--408.
\bibitem{chartier} Chartier, P., Faou, E., and Murua, M., An algebraic approach to invariant preserving integators: the case of quadratic and Hamiltonian invariants, {\em Numer. Math. \bf 103} (2006), 575--590. 

\bibitem{cmoq} Celledoni, E., McLachlan, R. I., Owren, B. and Quispel, G. R. W., Energy-preserving integrators and the structure of B-series, {\em  Foundations of Computational Mathematics \bf 10} (2010), 673--693.

\bibitem{co-kr} Connes, A. and Kreimer, D., Lessons from quantum field theory: Hopf algebras and spacetime geometries. {\em Letters in Mathematical Physics \bf 48} (1999), 85-96.

\bibitem{dur} A. D\"ur, {\em M\"obius functions, incidence algebras and power series representations}, Springer, Berlin 1986, pp. 88--90.

\bibitem{fauvet} Fauvet, F., Menous, F. and Sauzin, D., Explicit linearization of one-dimensional germs through tree-expansions, preprint, 2014. 
\bibitem{gerstenhaber63} Gerstenhaber, M., The cohomology structure of an associative ring, {\em Ann.  Math. \bf 78} (1963), 267--288.
\bibitem{gubinello} Gubinelli, M.,  Ramification of rough paths, {\em Journal of Differential Equations \bf 248}  (2010), 693--721.

\bibitem{hlw} Hairer, E., Lubich, C., and Wanner, G., {\em Geometric numerical integration: structure-preserving algorithms for ordinary differential equations}, 2nd ed., Springer, Berlin, 2006.
%\bibitem{ha-wa1} Hairer, E. and Wanner, G., Multistep-multistage-multiderivative methods for ordinary differential equations, {\em Computing \bf 11} (1973), 287--303.

\bibitem{ha-wa} Hairer, E., and Wanner, G., On the Butcher group and general multi-value methods, {\em	Computing \bf 13} (1974), 1--15.

\bibitem{ha-no-wa} Hairer, E., N\o rsett, S. P., \& Wanner, G., {\em Solving ordinary differential equation I: Nonstiff problems},  Springer, Berlin, 1987.
\bibitem{harary} Harary, F., Independent discoveries in graph theory, {\em Ann. New York Acad. Sci. \bf 328} (1979), 1--4.



\bibitem{iserles} Iserles, A., {\em A first course in the numerical analysis of differential equations}, Cambridge University Press, Cambridge, 2009.
\bibitem{is-qu-ts} Iserles, A., Quispel, G. R. W., and Tse, P. S. P., B-series methods cannot be volume-preserving, {\em BIT Numerical Mathematics \bf 47} (2007), 351-378.

\bibitem{jelliss} Jelliss, G. P., Knight's Tour notes, {\tt http://www.mayhematics.com/t/2n.htm}.
\bibitem{heleking} King-Hele, D., {\em  A Tapestry of Orbits}, Cambridge University Press, 2005.
\bibitem{kh} King-Hele, D., The destruction of the Royal Aircraft Establishment, \\ {\tt https://www.youtube.com/watch?v=E0fSLiAa9Zw}.

\bibitem{khashin} Khashin, S., Butcher algebras for Butcher systems, {\em Numer. Alg. \bf 63} (2013), 679--689.
\bibitem{mmmv} McLachlan, R. I., Modin, K., Munthe-Kaas, H., and Verdier, O., B--series are exactly the affine-equivariant methods, {\em Numer. Math.} (2015), pp. 1--24.
\bibitem{merson} Merson, R. H., An operational method for the study of integration processes, in {\em Proceedings
of Conference on Data Processing and Automatic Computing Machines vol. 1},  Weapons
Research Establishment, Salisbury, South Australia, 1957, pp. 1--25.
\bibitem{mk-lu} Munthe-Kaas, H. Z. and Lundervold, A., On post-Lie algebras, Lie--Butcher series and moving frames. {\em Found.\ Comp.\ Math. \bf 13} (2013), 583--613.
\bibitem{mko} Munthe-Kaas, H. and Verdier, O. (2015). Aromatic Butcher series. {\em Found.\ Comp.\ Math. \bf 16} (2016), 183--215.
\bibitem{sa-mu} Sanz-Serna, J. M. and Murua, A., Formal series and numerical integrators: some history and some new techniques, in {\em Proceedings of the 8th International Congress on Industrial and Applied Mathematics (ICIAM 2015)}, Lei Guo and Zhi-Ming eds., Higher Education Press, Beijing, 2015, 311--331.

\bibitem{vinberg} Vinberg, E. B., The theory of convex homogeneous cones, {\em Trans. Moscow Math. Soc. \bf 12} (1963), 340--403.
\bibitem{wanner} Wanner, G., personal communication.
\end{thebibliography}
\end{document}